\documentclass[11pt]{article}
\usepackage{amsmath}
\usepackage{pxfonts}
\usepackage{yfonts}
\usepackage{dsfont}
\usepackage{marvosym}
\usepackage{graphicx}
\usepackage{relsize}

\def\bel{\begin{equation}\label}
\def\eeq{\end{equation}}
\def\ds{\displaystyle}

\def\mt{\longrightarrow}
\def\v{\vskip 1em}

\def\ve{\varepsilon}
\def\R{\mathds R}
\def\Z{\mathds Z}
\def\C{\mathfrak{B}}

\def\Cx{\mathds C}
\def\N{{\bf N}}

\def\Re{{\bf Re}}
\def\Im {{\bf Im}}
\def\S{{\bf S}}

\def\O{{\bf O}}
\def\P{{\bf P}}

\def\J{{\bf J}}

\def\L{{\bf L}}

\def\p{{\partial}}

\def\i{{\bf i}}

\def\Hat{\widehat}

\def\I{{\bf I}}

\def\q{\mathfrak{q}}

\def\alpha{\alphaup}
\def\beta{\betaup}
\def\xi{{\xiup}}
\def\eta{{\etaup}}
\def\tau{{\tauup}}
\def\rho{{\rhoup}}
\def\phi{{\phiup}}
\def\psi{{\psiup}}
\def\omega{\omegaup}
\def\varphi{{\varphiup}}
\def\lambda{{\lambdaup}}
\def\c{{\bf c}}

\def\m{{\bf m}}
\def\n{{\bf n}}

\def\a{{\bf a}}
\def\b{{\bf b}}

\def\p{\partial}

\def\q{{\deltaup}}

\newtheorem{remark}{Remark}[section]

\begin{document}
 \[\begin{array}{cc}\hbox{\LARGE{\bf A family of convolution operators, part two}}
\end{array}\]

  \[\hbox{Zipeng Wang}\]

 \begin{abstract}
 We study a family of convolution operators. Their regarding Fourier multipliers are defined in terms of distributions having singularity on the light-cone in $\R^{n+1}$.
 \end{abstract}

 \section{Introduction}
 \setcounter{equation}{0} 
As investigated by Gelfand and Shilov \cite{Gelfand-Shilov},
 $\Lambda^\alpha, \alpha\in\Cx$ is a distribution defined in $\R^{n+1}$ by analytic continuation from
\bel{Lambda^alpha}
\begin{array}{cc}\ds
\hbox{\small{$\Re\alpha>{n-1\over2}$}},
\qquad
 \Lambda^\alpha(x,t)~=~ \hbox{\small{$\pi^{\alpha-{n+1\over2}}\Gamma^{-1}\left(\alpha-{n-1\over 2}\right)$}}
  \left({1\over t^2-|x|^2}\right)^{{n+1\over 2}-\alpha}_+. 
  \end{array}
\eeq
For $0<\Re\alpha<1$, 
the Fourier transform of $\Lambda^\alpha$ agrees on the function
\bel{Lambda Transform}
\begin{array}{rl}
  \Hat{\Lambda}^\alpha(\xi,\tau)~=~
\pi^{{n-1\over 2}-2\alpha}\Gamma\left(\alpha\right)
 \left\{ ~{\ds \left({1\over \tau^2-|\xi|^2}\right)^\alpha_-} - \sin\pi\left(\alpha-{1\over 2}\right){\ds \left( {1\over \tau^2- |\xi|^2}\right)^\alpha_+}~\right\}
  \end{array}
 \eeq
whenever $|\tau|\neq|\xi|$. 
See appendix B. 

Define
 \bel{I alpha}
 \begin{array}{lr}\ds
\I^\alpha f(x)
~=~\int_{\R^n} e^{2\pi\i x\cdot \xi} \Hat{f}(\xi) \Hat{\phi}(\xi)\left\{\int_0^1 \Hat{\Lambda}^\alpha(\xi,\tau) \tau d\tau\right\}d\xi
\end{array}
\eeq
where $\Hat{\phi}$ is a radially symmetric, smooth cut-off function supported in the shell: ${1\over 3}<|\xi|\leq3$.

$\diamond$  {\small Throughout,  $\C>0$ is regarded as a generic constant whose value depends on its sub-indices. }  

$\diamond$ {\small We write $\c>0$ for some fixed, large enough constant.}
\v
{\bf Theorem One} ~~{\it For ${1\over 2}<\Re\alpha<1$,   we have
 \bel{Result Two} 
 \begin{array}{cc}\ds
  \left\| \I^\alpha f\right\|_{\L^p(\R^n)} ~\leq~\C_{\Re\alpha}~e^{\c|\Im\alpha|}~\left\| f\right\|_{\L^p(\R^n)}, \qquad  {n-1\over 2n}~\leq ~{1\over p}~\leq~ {n+1\over 2n}.  
 \end{array}
  \eeq
  }\\
 This paper is organized as follows. In the next section, we show that $\I^\alpha$ defined in (\ref{I alpha}) is the same operator introduced in our previous work \cite{Wang}. In section 3, we give an implication of {\bf Theorem One} to the $\L^p$-boundedness of Bochner-Riesz operator.
 We add 2 appendices in the end for certain preliminary estimates.

 \section{Formulation on the main result}
 \setcounter{equation}{0}
 As  investigated by Strichartz \cite{Strichartz}, 
$\Omega^\alpha$ is a distribution defined  by analytic continuation from
\bel{Omega^alpha}
\begin{array}{cc}\ds
\hbox{\small{$\Re\alpha>{n-1\over 2}$}},\qquad \Omega^{\alpha}(x)~=~\hbox{\small{$\pi^{\alpha-{n+1\over 2}}\Gamma^{-1}\left(\alpha-{n-1\over 2}\right)$}} 
 \left({1 \over1-|x|^2}\right)^{{n+1\over 2}-\alpha}_+.
 \end{array}
\eeq
Equivalently,  it can be  defined by
 \bel{Omega^alpha Transform} 
\begin{array}{lr}\ds
\Hat{\Omega}^{\alpha}(\xi)~=~
\left({1\over|\xi|}\right)^{{n\over 2}-\big[{n+1\over 2}-\alpha\big]} \J_{{n\over 2}-\big[{n+1\over 2}-\alpha\big]}\Big(2\pi|\xi|\Big)
\\\\ \ds~~~~~~~~~
~=~
\left({1\over|\xi|}\right)^{\alpha-{1\over 2}} \J_{\alpha-{1\over 2}}\Big(2\pi|\xi|\Big),\qquad \alpha\in\Cx.
\end{array}
\eeq
See (\ref{Omega^z Transform}) in appendix A.
By using the integral formula of Bessel functions in (\ref{Bessel}), we  have
 \bel{r local inte}
 \Hat{\Omega}^\alpha (\xi)~=~ \pi^{\alpha-1}\Gamma\left(\alpha\right)\int_{-1}^1 e^{2\pi\i |\xi|s } (1-s^2)^{\alpha-1} ds.
  \eeq 
Recall $\Lambda^\alpha, \alpha\in\Cx$ which is a distribution defined in $\R^{n+1}$ by analytic continuation from (\ref{Lambda^alpha}). 
We find
\bel{Kernel singular}
 \begin{array}{lr}\ds
 h\ast\Lambda^\alpha(x,t)~=~\hbox{\small{$\pi^{\alpha-{n+1\over2}}\Gamma^{-1}\left(\alpha-{n-1\over 2}\right)$}}\iint_{|u|<|r|} h(x-u,t-r) 
  \left({1\over r^2-|u|^2}\right)^{{n+1\over 2}-\alpha}dudr 
 \\\\ \ds~~~~~~~~~~~~~~~~~~
 ~=~\iint_{\R^{n+1}} h(x-u,t-r)\Omega^\alpha\left({u\over r}\right)|r|^{2\alpha-1-n}dudr
  \\\\ \ds~~~~~~~~~~~~~~~~~~
 ~=~\lim_{N\mt\infty}\iint_{\R^{n+1}} e^{2\pi\i \big[x\cdot\xi+ t\tau\big]}\Hat{h}(\xi,\tau) \left\{ \int_{-N}^N e^{-2\pi\i \tau r}\Hat{\Omega}^\alpha (r\xi) |r|^{2\alpha-1} dr\right\} d\xi d\tau 
 \end{array}
\eeq
whenever $\Re\alpha>{n-1\over 2}$.
 \begin{remark}
By the principle of analytic continuation, we must have  (\ref{Kernel singular}) hold for every $\Re\alpha>0$.
\end{remark} 
 Let $0<\Re\alpha<1$. By using the asymptotic expansion of Bessel functions in (\ref{J asymptotic})-(\ref{J O}), write
\bel{Lambda transform Rewrite}
\begin{array}{lr}\ds
\int_\R e^{-2\pi\i \tau r}\Hat{\Omega}^\alpha (r\xi) |r|^{2\alpha-1} dr
~=~\int_\R e^{-2\pi\i \tau r}\left({1\over|r\xi|}\right)^{\alpha-{1\over 2}} \J_{\alpha-{1\over 2}}\Big(2\pi|r\xi|\Big) |r|^{2\alpha-1} dr\qquad\hbox{\small{by (\ref{Omega^alpha Transform})}}
\\\\ \ds~~~~~~~~~~~~~~~~~~~~~~~~~~~~~~~~~~~~~~~~~~
~=~{1\over \pi}\int_\R e^{-2\pi\i \tau r}\left({1\over|r\xi|}\right)^\alpha \cos\left[2\pi|r\xi|-{\pi\over 2}\alpha\right] |r|^{2\alpha-1} dr
\\\\ \ds~~~~~~~~~~~~~~~~~~~~~~~~~~~~~~~~~~~~~~~~~~
~+~\int_\R e^{-2\pi\i \tau r}\mathcal{E}^\alpha\left(2\pi |r\xi|\right) |r|^{2\alpha-1} dr

\\\\ \ds~~~~~~~~~~~~~~~~~~~~~~~
\hbox{\small{where}} \qquad\left|\mathcal{E}^\alpha(\rho)\right|~\leq~\C_{\Re\alpha}~e^{\c|\Im\alpha|}\left\{\begin{array}{lr}\ds \rho^{-\Re\alpha},\qquad 0<\rho\leq1,
 \\ \ds
 \rho^{-\Re\alpha-1},\qquad \rho>1.
 \end{array}\right.
 \end{array}
 \eeq 
 Consider $1/10\leq|\xi|\leq10$. Observe that the last integral in (\ref{Lambda transform Rewrite}) has  norm bounded by $\C_{\Re\alpha}~e^{\c|\Im\alpha|}$.
Because of Euler's formulae, we replace the cosine function in (\ref{Lambda transform Rewrite}) with $e^{-\i\left({\pi\over 2}\right)\alpha} e^{2\pi\i |r\xi|}$ or $e^{\i\left({\pi\over 2}\right)\alpha} e^{-2\pi\i |r\xi|}$.  

By  integration by parts $w.r.t~r$, we find
 \bel{main term norm r}
 \begin{array}{lr}\ds
 \left| \int_{2^{j-1}}^{2^j} e^{2\pi\i \big[|\xi|-\tau\big] r}\left({1\over|r\xi|}\right)^\alpha  |r|^{2\alpha-1} dr\right|
 \\\\ \ds
~\leq~ \C_{\Re\alpha}~e^{\c|\Im\alpha|}  \left\{\begin{array}{lr} 2^{j\Re\alpha} ,~~~~~~~~~~~~~~~~~~~~~~\qquad 2^j\leq \left||\tau|-|\xi|\right|^{-1},
 \\ \ds
 \left|{1\over |\tau|-|\xi|}\right| 2^{j\big[\Re\alpha-1\big]},\qquad 2^j> \left||\tau|-|\xi|\right|^{-1}. 
 \end{array}\right.
 \end{array}
 \eeq
 By symmetry reason,  (\ref{main term norm r}) implies
 \bel{Int Norm Est r}
 \begin{array}{lr}\ds
 \left| \int_\R e^{-2\pi\i \tau r}\left({1\over|r\xi|}\right)^\alpha \cos\left[2\pi|r\xi|-{\pi\over 2}\alpha\right] |r|^{2\alpha-1} dr\right|
 \\\\ \ds
 ~\leq~ \sum_{j\in\Z}\left| \int_{2^{j-1}\leq|r|<2^j} e^{-2\pi\i \tau r}\left({1\over|r\xi|}\right)^\alpha \cos\left[2\pi|r\xi|-{\pi\over 2}\alpha\right] |r|^{2\alpha-1} dr\right| 
 \\\\ \ds
 ~\leq~  \C_{\Re\alpha}~e^{\c|\Im\alpha|} \sum_{2^j\leq \left||\tau|-|\xi|\right|^{-1}} 2^{j\Re\alpha}+\C_{\Re\alpha}~e^{\c|\Im\alpha|} \sum_{2^j> \left||\tau|-|\xi|\right|^{-1}}  \left|{1\over |\tau|-|\xi|}\right| 2^{j\big[\Re\alpha-1\big]}
 \\\\ \ds
  ~\leq~  \C_{\Re\alpha}~e^{\c|\Im\alpha|}~ \left|{1\over |\tau|-|\xi|}\right|^{\Re\alpha}.
  \end{array}  
 \eeq
On the other hand, the Fourier transform of $\Lambda^\alpha, 0<\Re\alpha<1$ agrees with the function
 $\Hat{\Lambda}^\alpha(\xi,\tau)$ in (\ref{Lambda Transform}) whenever $|\tau|\neq|\xi|$.  
 From (\ref{Kernel singular}) to (\ref{Int Norm Est r}), by applying Lebesgue's dominated convergence theorem, we have
  \bel{Kernel singular region}
 \begin{array}{lr}\ds
 h\ast\Lambda^\alpha(0,0)
 ~=~\iint_{\R^{n+1}} \Hat{h}(\xi,\tau) 
 \left\{ \int_\R e^{-2\pi\i \tau r}\Hat{\Omega}^\alpha (r\xi) |r|^{2\alpha-1} dr\right\} d\xi d\tau 
 \\\\ \ds~~~~~~~~~~~~~~~~~~
 ~=~\iint_{\R^{n+1}} \Hat{h}(\xi,\tau) 
\Hat{\Lambda}^\alpha(\xi,\tau) d\xi d\tau 
 \end{array}
\eeq
 for every Schwartz function $h$ of which $\Hat{h}(\xi,\tau)$ is supported in the cylinder: ${1\over 10}\leq|\xi|\leq10$.

Let $\varphi\in\mathcal{C}^\infty_o(\R)$ be a smooth bump-function  such that $\varphi(t)=1$ if $|t|\leq1$ and $\varphi(t)=0$ if $|t|>2$.  Given $(\eta,s)\in\R^n\times\R$ for $|s|\neq|\eta|$ and ${1\over 3}\leq|\eta|\leq3$, we consider a family of {\it good kernels}: $\Hat{h}_\ve(\xi,\tau)=\c_\varphi  \ve^{-(n+1)}\varphi \left[\ve^{-1}\sqrt{|\xi-\eta|^2+(\tau-s)^2}\right],~0<\ve<{1\over 10}$
  where $\ds\c_\varphi^{-1}=\int_\R \varphi(t)dt$.  
  
Replace $\Hat{h}(\xi,\tau)$ by $\Hat{h}_\ve(\xi,\tau)$ in (\ref{Kernel singular region}). By taking $\ve\mt0$, we find
\bel{Lambda transform rewrite}
\begin{array}{lr}\ds
\Hat{\Lambda}^\alpha(\xi,\tau)~=~\int_\R e^{-2\pi\i \tau r}\Hat{\Omega}^\alpha (r\xi) |r|^{2\alpha-1} dr
\end{array} 
 \eeq
 for $|\tau|\neq|\xi|$ and ${1\over 3}\leq|\xi|\leq3$.

Recall $\I^\alpha f$ defined in (\ref{I alpha}). By applying Lebesgue's dominated convergence theorem again, we have
 \bel{I alpha rewrite}
 \begin{array}{cc}\ds
 \Hat{\I^\alpha f}(\xi)~=~ \Hat{f}(\xi)\Hat{\P}^\alpha(\xi),
 \\\\ \ds
  \Hat{\P}^\alpha(\xi)~=~\Hat{\phi}(\xi)\int_0^1   \Hat{\Lambda}^\alpha(\xi,\tau)\tau d\tau  
  ~~~~~~~~~~~~~~~~~~~~~~~~~~~~~~~~~~~~~~~~~
  \\\\ \ds~~~~~~~~~~~~~~~~~~~~~~~~~~~~~~~~~~
  ~=~ \Hat{\phi}(\xi)\int_{0<\tau<1,~\tau\neq|\xi|}  \left\{\int_\R e^{-2\pi\i \tau r}\Hat{\Omega}^\alpha (r\xi) |r|^{2\alpha-1} dr\right\}\tau d\tau \qquad\hbox{\small{by (\ref{Lambda transform rewrite})}}  
    \\\\ \ds~~~~~~~~~~~~~~~~~~~~~~~~
  ~=~\Hat{\phi}(\xi)\lim_{N\mt\infty}\int_{0<\tau<1,~\tau\neq|\xi|} \left\{ \int_{-N}^N e^{-2\pi\i \tau r}\Hat{\Omega}^\alpha (r\xi) |r|^{2\alpha-1} dr\right\} \tau d\tau
  \\\\ \ds~~~~
~=~\Hat{\phi}(\xi) \int_\R \left\{\int_0^1 e^{-2\pi\i \tau r} \tau d\tau\right\} \Hat{\Omega}^\alpha (r\xi) |r|^{2\alpha-1} dr.  ~~~~~
 \end{array}
 \eeq
Let
 \bel{omega}
 \begin{array}{lr}\ds
 \omega(r)~=~e^{2\pi\i r}\int_0^1 e^{-2\pi\i \tau r} \tau d\tau
 \\\\ \ds~~~~~~~~
 ~=~{-1\over 2\pi\i}{1\over r}-{1\over 4\pi^2r^2} \left[1-e^{-2\pi\i r}\right].
 \end{array}
 \eeq 
  \begin{remark} We find $|\omega(r)|\leq\C\Big[1+|r|\Big]^{-1}$.
 \end{remark}
 From (\ref{I alpha rewrite})-(\ref{omega}), we write
\bel{P rewrite} 
\begin{array}{cc}\ds
\I^\alpha f(x)~=~\int_{\R^n} e^{2\pi\i x\cdot\xi} \Hat{f}(\xi)\Hat{\P}^\alpha(\xi)d\xi,
\\\\ \ds
\Hat{\P}^\alpha(\xi)~=~\Hat{\phi}(\xi)\int_\R e^{-2\pi\i  r}\Hat{\Omega}^\alpha (r\xi)\omega(r) |r|^{2\alpha-1} dr.
\end{array}
\eeq 
\begin{remark}
The proof of {\bf Theorem One} for $\I^\alpha f, {1\over 2}<\Re\alpha<1$ defined as 
(\ref{P rewrite}) is given in the paper \cite{Wang}.
 \end{remark}

 \section{Implication to the Bochner-Riesz summability}
 \setcounter{equation}{0}
  The Bochner-Riesz operator is defined by
   \bel{S}
  \S^\q f(x)~=~\int_{\R^n}e^{2\pi\i x\cdot\xi} \Hat{f}(\xi) \left( 1-|\xi|^2\right)^\q_+ d\xi,\qquad \Re\q\ge0.
  \eeq
At $\q=0$, we revisit on the famous ball multiplier problem.  $\S^0$ is bounded only on $\L^2(\R^n)$ proved by Fefferman \cite{Fefferman}.

 Let $0<\Re\q<{1\over 2}$. 
Consider
  \bel{S_phi}
  \begin{array}{lr}\ds
  \S^\q_\psi f(x)~=~\int_{\R^n} e^{2\pi\i x\cdot\xi} \Hat{f}(\xi) \Hat{\psi}(\xi)\left(1-|\xi|^2\right)^\q_+ d\xi
 \end{array}  
  \eeq  
where $ \Hat{\psi}(\xi)=\left[\Hat{\phi}(\xi)\right]^2$ as shown in (\ref{I alpha}).

From direct computation, we find
 \bel{Marcinkiewicz multiplier}
 \left|\left(\xi\p_\xi\right)^\gamma\left[1-\Hat{\psi}(\xi)\right] \left(1-|\xi|^2\right)^\q_+\right|~\leq~\C_{\Re\q~\gamma}~e^{\c|\Im\q|}
 \eeq
 for every multi-index $\gamma$. Consequently, $\left[1-\Hat{\psi}(\xi)\right] \left(1-|\xi|^2\right)^\q_+$
is a $\L^p$-Fourier multiplier for $1<p<\infty$.  
 \v
  {\bf Theorem Two}~~{\it Let $\S^\q_\psi$ defined in (\ref{S_phi}) for $0<\Re\q<{1\over 2}$. We have
 \bel{RESULT ONE}
 \left\| \S^\q_\psi f\right\|_{\L^p(\R^n)}~\leq~\C_{\Re\q~p}~e^{\c|\Im\q|}~\left\| f\right\|_{\L^p(\R^n)},\qquad  {n-1\over 2n}~\leq~{1\over p}~\leq~{n+1\over 2n}. 
 \eeq} 
 
 Next, we show that {\bf Theorem one} implies {\bf Theorem Two}. 
 
 Let $1-\q=\alpha$. We define
\bel{m+}
\begin{array}{lr}\ds
\m^\alpha_+(\xi)~=~\Hat{\phi}(\xi)\int_0^1  \left({1\over \tau^2-|\xi|^2}\right)^\alpha_+ \tau d\tau
\\\\ \ds~~~~~~~~~~
~=~\Hat{\phi}(\xi)\left\{ \begin{array}{lr}\ds
\int_{|\xi|}^1 \left({1\over \tau^2-|\xi|^2}\right)^\alpha \tau d\tau,\qquad |\xi|<1
\\ \ds~~~~~~~~~~~~~~~~~~
0\qquad~~~~~~~~~~~~~~~ |\xi|\ge1
\end{array}\right.
\\\\ \ds~~~~~~~~~~
~=~{1\over 2} (1-\alpha)^{-1} \Hat{\phi}(\xi) \left(1-|\xi|^2\right)_+^{1-\alpha}
  \end{array}
  \eeq
and
\bel{m-} 
\begin{array}{lr}\ds
\m^\alpha_-(\xi)~=~\Hat{\phi}(\xi)\int_0^1  \left({1\over \tau^2-|\xi|^2}\right)^\alpha_- \tau d\tau 
\\\\ \ds~~~~~~~~~~
~=~(-1)^{-\alpha}\Hat{\phi}(\xi)\left\{ \begin{array}{lr}\ds
\int_0^{|\xi|} \left({1\over \tau^2-|\xi|^2}\right)^\alpha\tau d\tau,\qquad |\xi|\leq1
\\ \ds
\int_0^1 \left({1\over \tau^2-|\xi|^2}\right)^\alpha\tau d\tau,\qquad~ |\xi|>1
\end{array}\right.
\\\\ \ds~~~~~~~~~~
 ~=~{1\over 2}(1-\alpha)^{-1}  \Hat{\phi}(\xi) \left[|\xi|^{2(1-\alpha)}-\left(1-|\xi|^2\right)_-^{1-\alpha}\right].
\end{array}
\eeq
Recall $\Hat{\Lambda}^\alpha(\xi,\tau) $ given at (\ref{Lambda Transform}).
From (\ref{m+})-(\ref{m-}), we find
 \bel{Lambda m}
 \begin{array}{lr}\ds
  \Hat{\phi}(\xi)\int_0^1  \Hat{\Lambda}^\alpha(\xi,\tau) \tau d\tau
  ~=~\pi^{{n-1\over 2}-2\alpha}\Gamma(\alpha) \Bigg\{ \m^\alpha_-(\xi)-\hbox{$\sin\pi\left(\alpha-{1\over 2}\right)$}\m^\alpha_+(\xi)\Bigg\}.
  \end{array}
 \eeq 
Moreover, define
 \bel{m^alpha}
 \m^\alpha(\xi)~=~\Hat{\phi}(\xi) \Bigg\{  -\left(1-|\xi|^2\right)_-^{1-\alpha}-\hbox{$\sin\pi\left(\alpha-{1\over 2}\right)$} \left(1-|\xi|^2\right)_+^{1-\alpha}\Bigg\}.
 \eeq 
 From (\ref{m+})-(\ref{m-}), (\ref{Lambda m}) and (\ref{m^alpha}), we have
  \bel{m alpha}
 \begin{array}{lr}\ds
  \Hat{\phi}(\xi)\int_0^1  \Hat{\Lambda}^\alpha(\xi,\tau) \tau d\tau
  ~=~
  \pi^{{n-2\over 2}-2\alpha}{1\over 2}(1-\alpha)^{-1} \Gamma(\alpha)\Big[  \m^\alpha(\xi)+\Hat{\phi}(\xi)  |\xi|^{2(1-\alpha)} \Big].
  \end{array}
 \eeq 
   \begin{remark}
$\m^z(\xi), z\in\Cx$ is called a restricted $\L^p$-Fourier multiplier if the regarding convolution is bounded on $\L^p(\R^n)$ for ${n-1\over 2n}\leq{1\over p}\leq{n+1\over 2n}$ with an implied constant $\C_{\Re z}e^{\c|\Im z|}$.
 \end{remark}  
Suppose ${1\over 2}<\Re\alpha<1$. Recall $\I^\alpha$ defined in (\ref{I alpha}) of which $\I^\alpha f$ satisfies  the $\L^p$-norm inequality in (\ref{Result Two}).  Because $\Hat{\phi}(\xi)|\xi|^{2(1-\alpha)}$ inside (\ref{m alpha}) is a $\L^p$-Fourier multiplier for $1<p<\infty$, we find $\m^\alpha(\xi)$  defined in (\ref{m^alpha}) as a restricted $\L^p$-Fourier multiplier.

Consider
 \bel{m_1}
\Hat{\phi}(\xi)\m^\alpha(\xi)~=~ \Hat{\psi}(\xi) \Bigg\{  -\left(1-|\xi|^2\right)_-^{1-\alpha}-\hbox{$\sin\pi\left(\alpha-{1\over 2}\right)$} \left(1-|\xi|^2\right)_+^{1-\alpha}\Bigg\}
\eeq
and
\bel{m_2}
\m^{{1\over 2}+{\alpha\over 2}}(\xi)~=~\Hat{\phi}(\xi) \Bigg\{  -\left(1-|\xi|^2\right)_-^{{1\over 2}-{\alpha\over 2}}-\hbox{$\sin\pi\left({\alpha\over 2}\right)$} \left(1-|\xi|^2\right)_+^{{1\over 2}-{\alpha\over 2}}\Bigg\}.
\eeq
  Clearly, 
 both $\Hat{\phi}(\xi)\m^\alpha(\xi)$ and $\m^{{1\over 2}+{\alpha\over 2}}(\xi)$ 
 are two  restricted $\L^p$-Fourier multipliers. 
 
 Furthermore, a direct computation shows
\bel{m_3}
\begin{array}{lr}\ds
\left[\m^{{1\over 2}+{\alpha\over 2}}(\xi)\right]^2
~=~ \Hat{\psi}(\xi)  \Bigg\{  \left(1-|\xi|^2\right)_-^{1-\alpha}+\hbox{$\sin^2\pi\left({\alpha\over 2}\right)$}\left(1-|\xi|^2\right)_+^{1-\alpha}\Bigg\}  
 \end{array}
\eeq 
which is indeed a restricted $\L^p$-Fourier multiplier.

By adding $\left[\m^{{1\over 2}+{\alpha\over 2}}(\xi)\right]^2$ and $\Hat{\phi}(\xi)\m^\alpha(\xi)$ together, we find
 \bel{subtract}
 \begin{array}{lr}\ds
\Hat{\psi}(\xi) \Big[\hbox{$\sin^2\pi\left({\alpha\over 2}\right)-\sin\pi\left(\alpha-{1\over 2}\right)$}\Big] \left(1-|\xi|^2\right)_+^{1-\alpha}
\end{array}
\eeq 
where $\sin^2\pi\left({\alpha\over 2}\right)-\sin\pi\left(\alpha-{1\over 2}\right)$ is non-zero for $0<\Re\alpha<1$.  This shows that
\bel{restricted multiplier} 
\Hat{\psi}(\xi)\left(1-|\xi|^2\right)_+^{1-\alpha}
\eeq
  is another restricted $\L^p$-Fourier multiplier.

\appendix 
  \section{Some estimates regarding Bessel functions}
 \setcounter{equation}{0}
$\bullet$ For  $a>-{1\over 2}, b\in\R$ and $\rho>0$,  a Bessel function $\J_{a+\i b}$ has an integral formula 
\bel{Bessel}
\J_{a+\i b}(\rho)~=~{(\rho/2)^{a+\i b}\over \pi^{1\over 2}\Gamma\left(a+{1\over 2}+\i b\right)} \int_{-1}^1 e^{\i\rho s } (1-s^2)^{a-{1\over 2}+\i b} ds.
\eeq
$\bullet$ For every $a>-{1\over 2}, b\in\R$ and $\rho>0$, we have 
\bel{J asymptotic}
\begin{array}{lr}\ds
\J_{a+\i b}(\rho) ~\sim~\left({\pi\rho\over 2}\right)^{-{1\over 2}} \cos\left[\rho-(a+\i b) {\pi\over 2}-{\pi\over 4}\right]
\\\\ \ds~~~~~~~~~~~~
~+~\left({\pi\rho\over 2}\right)^{-{1\over 2}} \sum_{k=1}^\infty \cos\left[\rho-(a+\i b) {\pi\over 2}-{\pi\over 4}\right] \a_k \rho^{-2k}+\sin\left[\rho-(a+\i b) {\pi\over 2}-{\pi\over 4}\right]~ \b_k \rho^{-2k+1},
\\\\ \ds
\a_k~=~(-1)^k[a+\i b, 2k]2^{-2k},\qquad \b_k~=~(-1)^{k+1}[a+\i b, 2k-1]2^{-2k+1},
\\\\ \ds
[a+\i b,m]~=~{\Gamma\left({1\over 2}+a+\i b+m\right)\over m!\Gamma\left({1\over 2}+a+\i b-m\right)},\qquad  m~=~0,1,2,\ldots
\end{array}
\eeq
in the sense of that
\bel{J O}
\begin{array}{lr}\ds
\left({d\over d\rho}\right)^\ell\left[ \J_{a+\i b}(\rho) ~-~\left({\pi\rho\over 2}\right)^{-{1\over 2}} \cos\left[\rho-(a+\i b) {\pi\over 2}-{\pi\over 4}\right]\right.
\\\\ \ds~~~
\left.~-~\left({\pi\rho\over 2}\right)^{-{1\over 2}}\sum_{k=1}^N  \cos\left[\rho-(a+\i b) {\pi\over 2}-{\pi\over 4}\right] \a_k \rho^{-2k}+\sin\left[\rho-(a+\i b) {\pi\over 2}-{\pi\over 4}\right]~ \b_k \rho^{-2k+1}\right]
\\\\ \ds
~=~\O\left(\rho^{-2N-{1\over 2}}\right)\qquad N\ge0,\qquad \ell\ge0,\qquad \rho\mt\infty
\end{array}
\eeq
where the implied constant is bounded by $\C_a e^{\c|b|}$.

$\bullet$  For every $a,b\in\R$ and $\rho>0$, 
\bel{J identity}
\J_{a-1+\i b} (\rho)~=~2\left[{a+\i b\over \rho}\right] \J_{a+\i b}(\rho)-\J_{a+1+\i b}(\rho).
\eeq
By using (\ref{Bessel}) and (\ref{J asymptotic}) together with (\ref{J identity}), we find the norm estimate in below.

$\bullet$ For  every $a, b\in\R$ and $\rho>0$, 
\bel{J norm}
\left| {1\over\rho^{a+\i b}}~\J_{a+\i b}(\rho)\right|~\leq~\C_a~\left({1\over 1+\rho}\right)^{{1\over 2}+a}~e^{\hbox{\small{\bf c}}|b|}.
\eeq  
More discussion of Bessel functions can be found in  the book of Watson \cite{Watson}.

Let $z\in\Cx$.  
$\Omega^z$  is a distribution defined by analytic continuation from
\bel{Omega^z}
\begin{array}{ccc}\ds
\Re z<1,\qquad\Omega^z(x)~=~\pi^{-z}\Gamma^{-1}\left(1-z\right)  \left({1\over 1-|x|^2}\right)^{z}_+.
\end{array}
\eeq

$\bullet$  $\Omega^z$ can be equivalently defined by 
 \bel{Omega^z Transform} 
\begin{array}{lr}\ds
\Hat{\Omega}^z(\xi)~=~\left({1\over|\xi|}\right)^{{n\over 2}-z} \J_{{n\over 2}-z}\Big(2\pi|\xi|\Big),\qquad z\in\Cx.
\end{array}
\eeq
From (\ref{Omega^z}), we have
\bel{Omega^lambda Fourier transform}
\begin{array}{lr}\ds
\Hat{\Omega}^z(\xi)~=~\pi^{-z}  \Gamma^{-1}\left(1-z\right)\int_{|x|<1} e^{-2\pi\i x\cdot\xi} \left({1\over 1-|x|^2}\right)^z dx
\\\\ \ds~~~~~~~~~
~=~\pi^{-z}  \Gamma^{-1}\left(1-z\right)\int_0^\pi\left\{\int_0^1 e^{-2\pi\i |\xi|r\cos\vartheta} (1-r^2)^{-z} r^{n-1}dr\right\} \omega_{n-2} \sin^{n-2}\vartheta d\vartheta
\\\\ \ds~~~~~~~~~
~=~\omega_{n-2}\pi^{-z}  \Gamma^{-1}\left(1-z\right)\int_{-1}^1\left\{\int_0^1 e^{2\pi\i|\xi|rs} (1-r^2)^{-z} r^{n-1}dr\right\} (1-s^2)^{n-3\over 2} ds~~~~ \hbox{\small{($-s=\cos\vartheta$)}}
\\\\ \ds~~~~~~~~~
~=~\omega_{n-2}\pi^{-z}  \Gamma^{-1}\left(1-z\right)\int_0^1\left\{\int_{-1}^1 e^{2\pi\i|\xi|rs} (1-s^2)^{n-3\over 2} ds\right\}(1-r^2)^{-z} r^{n-1}dr
\\\\ \ds~~~~~~~~~
~=~\omega_{n-2}\pi^{-z}  \Gamma^{-1}\left(1-z\right)\int_0^1\left\{\int_{-1}^1 \cos(2\pi|\xi|rs) (1-s^2)^{n-3\over 2} ds\right\}(1-r^2)^{-z} r^{n-1}dr
\end{array}
\eeq
where 
\bel{area}
\omega_{n-2}~=~2 \pi^{n-1\over 2} \Gamma^{-1}\left({n-1\over 2}\right).
\eeq
Recall the Beta function identity:
\bel{Beta}
{\Gamma(z)\Gamma(w)\over \Gamma(z+w)}~=~\int_0^1 r^{z-1}(1-r)^{w-1}dr,\qquad \Re z>0,~~~\Re w>0.
\eeq
By writing out the Taylor expansion of the cosine function inside (\ref{Omega^lambda Fourier transform}), we find
\bel{Omega^lambda Fourier transform Sum}
\begin{array}{lr}\ds
\omega_{n-2}\pi^{-z}  \Gamma^{-1}\left(1-z\right)\int_0^1\left\{\int_{-1}^1 \cos(2\pi|\xi|rs) (1-s^2)^{n-3\over 2} ds\right\}(1-r^2)^{-z} r^{n-1}dr
\\\\ \ds
=~\omega_{n-2}\pi^{-z}  \Gamma^{-1}\left(1-z\right)
\\ \ds~~~~~
\sum_{k=0}^\infty (-1)^k {(2\pi|\xi|)^{2k}\over (2k)!}\left\{\int_{-1}^1 s^{2k} (1-s^2)^{n-3\over 2} ds\right\}\left\{\int_0^1 r^{2k+n-1}(1-r^2)^{-z} dr\right\}
\\\\ \ds
=~{1\over 2}\omega_{n-2}\pi^{-z}  \Gamma^{-1}\left(1-z\right)
\\ \ds~~~~~
\sum_{k=0}^\infty (-1)^k {(2\pi|\xi|)^{2k}\over (2k)!}\left\{\int_0^1 t^{k+{1\over 2}-1} (1-t)^{{n-1\over 2}-1} dt\right\}\left\{\int_0^1\rho^{k+{n\over 2}-1}(1-\rho)^{1-z-1} d\rho\right\}
~~~~\hbox{\small{($t=s^2$,~$\rho=r^2$)}}
\\\\ \ds
=~{1\over 2}\omega_{n-2}\pi^{-z}  \Gamma^{-1}\left(1-z\right)\sum_{k=0}^\infty (-1)^k {(2\pi|\xi|)^{2k}\over (2k)!}~{\Gamma\left(k+{1\over 2}\right)\Gamma\left({n-1\over 2}\right)\over\Gamma\left(k+{n\over 2}\right)}~{ \Gamma\left(k+{n\over 2}\right) \Gamma\left(1-z\right)         \over    \Gamma\left(k+{n\over 2}+1-z\right)} ~~~~
 \hbox{\small{by (\ref{Beta})}}
\\\\ \ds
=~ \pi^{{n-1\over 2}-z}\sum_{k=0}^\infty (-1)^k {(2\pi|\xi|)^{2k}\over (2k)!}~{\Gamma\left(k+{1\over 2}\right)        \over    \Gamma\left(k+{n\over 2}+1-z\right)},\qquad\hbox{\small{by (\ref{area})}}.
\end{array}
\eeq
On the other hand, we have
\bel{Omega cos}
\begin{array}{lr}\ds
\left({1\over|\xi|}\right)^{{n\over 2}-z} \J_{{n\over 2}-z}\Big(2\pi|\xi|\Big)
~=~\hbox{\small{$\pi^{{n-1\over 2}-z}  \Gamma^{-1}\left({n+1\over 2}-z\right)$}}\int_{-1}^1 e^{2\pi\i  |\xi| s} (1-s^2)^{{n-1\over 2}-z} ds\qquad\hbox{\small{by (\ref{Bessel})}}
\\\\ \ds~~~~~~~~~~~~~
~=~\hbox{\small{$\pi^{{n-1\over 2}-z}  \Gamma^{-1}\left({n+1\over 2}-z\right)$}}\int_{-1}^1 \cos\left(2\pi  |\xi| s\right) (1-s^2)^{{n-1\over 2}-z} ds
\\\\ \ds~~~~~~~~~~~~~
~=~\hbox{\small{$\pi^{{n-1\over 2}-z}  \Gamma^{-1}\left({n+1\over 2}-z\right)$}} \sum_{k=0}^\infty (-1)^k {(2\pi|\xi|)^{2k}\over (2k)!}\int_{-1}^1 s^{2k} (1-s^2)^{{n-1\over 2}-z} ds
\\\\ \ds~~~~~~~~~~~~~
~=~\hbox{\small{$\pi^{{n-1\over 2}-z}  \Gamma^{-1}\left({n+1\over 2}-z\right)$}} \sum_{k=0}^\infty (-1)^k {(2\pi|\xi|)^{2k}\over (2k)!}\int_0^1 \rho^{k+{1\over 2}-1} (1-\rho)^{{n+1\over 2}-z-1} d\rho\qquad \hbox{\small{($\rho=s^2$)}}
\\\\ \ds~~~~~~~~~~~~~
~=~\pi^{{n-1\over 2}-z} \sum_{k=0}^\infty (-1)^k {(2\pi|\xi|)^{2k}\over (2k)!}~{\Gamma\left(k+{1\over 2}\right)    \over    \Gamma\left(k+{n\over 2}+1-z\right)}\qquad\hbox{\small{by (\ref{Beta}).}}
\end{array}
\eeq

\section{Fourier transform of $\Lambda^\alpha$}
\setcounter{equation}{0}
We derive the formula of $\Hat{\Lambda}^\alpha(\xi,\tau)$ in (\ref{Lambda Transform})
by following the lines in p. $253-284$,
Chapter III of Gelfand and Shilov \cite{Gelfand-Shilov}.

$\Hat{\hbox{U}}^\alpha$, $\Hat{\hbox{V}}^\alpha$ are distributions defined by analytic continuation from
\bel{UV}
\begin{array}{lr}\ds
\Hat{\hbox{U}}^\alpha(\xi,\tau)~=~\Gamma^{-1}\left(1-\alpha\right)
  \left({1\over \tau^2-|\xi|^2}\right)^\alpha_-,
\qquad
 \Hat{\hbox{V}}^\alpha(\xi,\tau)~=~\Gamma^{-1}\left(1-\alpha\right)
  \left( {1\over \tau^2-|\xi|^2}\right)^\alpha_+,
\qquad
  \hbox{\small{$\Re\alpha<1$}}.  
\end{array}
\eeq
For brevity, we denote  $\sigma(\alpha)={n+1\over 2}-\alpha,~ \alpha\in\Cx$. $\Pi^\alpha$, $\Lambda^\alpha$ are distributions defined by analytic continuation from
\bel{PiLambda}
\begin{array}{lr}\ds
\Pi^\alpha(x,t)~=~\Gamma^{-1}\left(1-\sigma(\alpha)\right)
  \left({1\over t^2-|x|^2}\right)^{\sigma(\alpha)}_-,
\qquad
 \Lambda^\alpha(x,t)~=~\Gamma^{-1}\left(1-\sigma(\alpha)\right)
  \left({1\over t^2-|x|^2}\right)^{\sigma(\alpha)}_+,
\\\\\ds~~~~~~~~~~~~~~~~~~~~~~~~~~~~~~~~~~~~~~~~~~~~~~~~~~~~~~~~~~~~~~~~~~~~~~~~~~~~~~~~~~~~~~~~~~~~~~~~~~~~~~~~~~~~~~~~~~~~~~~~~~~~
 \hbox{\small{$\Re\sigma(\alpha)<1$}}.  
  \end{array}
\eeq 
Let $a>0, b>0$. Define
\bel{rho xi,tau}
\rho(\xi,\tau,a,b)~=~\Bigg\{ \left[\tau^2-|\xi|^2\right]^2+\left[a \tau^2+b|\xi|^2\right]^2\Bigg\}^{1\over 2},\qquad \cos\theta(\xi,\tau,a,b)~=~{\tau^2-|\xi|^2\over\rho(\xi,\tau,a,b)}.
\eeq
For $\Re\alpha<{n+1\over 2}$, we assert
 \bel{P,R,ab}
 \begin{array}{lr}\ds
\Hat{\hbox{P}}^{\alpha~a~b}(\xi,\tau)~=~\left[ \tau^2-|\xi|^2+\i a \tau^2+\i b |\xi|^2\right]^{-\alpha}
\\\\ \ds~~~~~~~~~~~~~~~~~~
~=~\rho(\xi,\tau,a,b)^{-\alpha} e^{-\i\alpha \theta(\xi,\tau,a,b)},
\end{array}
~~~~
\begin{array}{lr}\ds
\Hat{\hbox{R}}^{\alpha~a~b}(\xi,\tau)~=~\left[ \tau^2-|\xi|^2-\i a \tau^2-\i b |\xi|^2\right]^{-\alpha}
\\\\ \ds~~~~~~~~~~~~~~~~~~
~=~\rho(\xi,\tau,a,b)^{-\alpha} e^{\i\alpha \theta(\xi,\tau,a,b)}.
\end{array}
\eeq
$\Hat{\hbox{P}}^{\alpha~a~b}$ and $\Hat{\hbox{R}}^{\alpha~a~b}$ are distributions defined by analytic continuation from (\ref{P,R,ab}). 
Recall (\ref{UV}). We find
\bel{PR}
\begin{array}{lr}\ds
 \Hat{\hbox{P}}^\alpha~\doteq~\lim_{a\mt0,~b\mt0} \Hat{\hbox{P}}^{\alpha~a~b}
 ~=~\Gamma\left(1-\alpha\right)\left[e^{-\i\pi \alpha}\Hat{\hbox{U}}^\alpha+ \Hat{\hbox{V}}^\alpha\right],
\\\\ \ds
\Hat{\hbox{R}}^\alpha~\doteq~\lim_{a\mt0,~b\mt0} \Hat{\hbox{R}}^{\alpha~a~b}~=~\Gamma\left(1-\alpha\right)\left[e^{\i\pi\alpha}\Hat{\hbox{U}}^\alpha+ \Hat{\hbox{V}}^\alpha\right].
\end{array}
\eeq
Let $\rho(x,t,a,b)$ and $\theta(x,t,a,b)$ defined as same as (\ref{rho xi,tau}).
For $\Re\sigma(\alpha)<{n+1\over 2}$, we consider
\bel{Phi,Psi,ab}
 \begin{array}{lr}\ds
\Phi^{\alpha~a~b}(x,t)~=~\left[ t^2-|x|^2+\i a t^2+\i b |x|^2\right]^{-\sigma(\alpha)}
\\\\ \ds~~~~~~~~~~~~~~~~~~
~=~\rho(x,t,a,b)^{-\sigma(\alpha)} e^{-\i\sigma(\alpha) \theta(x,t,a,b)},
\end{array}
~~~~
\begin{array}{lr}\ds
\Psi^{\alpha~a~b}(x,t)~=~\left[ t^2-|x|^2-\i a t^2-\i b |x|^2\right]^{-\sigma(\alpha)}
\\\\ \ds~~~~~~~~~~~~~~~~~~
~=~\rho(x,t,a,b)^{-\sigma(\alpha)} e^{\i\sigma(\alpha) \theta(x,t,a,b)}.
\end{array}
\eeq
$\Phi^{\alpha~a~b}$ and $\Psi^{\alpha~a~b}$ are distributions defined by analytic continuation from (\ref{Phi,Psi,ab}). 
Recall (\ref{PiLambda}). 
We have
\bel{PhiPsi} 
\begin{array}{lr}\ds
\Phi^\alpha~\doteq~\lim_{a\mt0,~b\mt0} \Phi^{\alpha~a~b}
~=~\Gamma\left(1-\sigma(\alpha)\right)\left[e^{-\i\pi \sigma(\alpha)} \Pi^\alpha+\Lambda^\alpha\right],
\\\\ \ds
\Psi^\alpha~\doteq~\lim_{a\mt0,~b\mt0} \Psi^{\alpha~a~b}
~=~\Gamma\left(1-\sigma(\alpha)\right)\left[e^{\i\pi \sigma(\alpha)} \Pi^\alpha+\Lambda^\alpha\right].
\end{array}
\eeq
$\bullet$ $\hbox{P}^\alpha$, $\hbox{R}^\alpha$ are analytic for $\alpha\in\Cx$ except the simple poles at 
$\alpha={n+1\over 2}+k,~ k=0,1,2,\ldots$.
  
$\bullet$ $\Phi^\alpha$, $\Psi^\alpha$ are analytic for $\alpha\in\Cx$ except the simple poles at 
$\lambda(\alpha)={n+1\over 2}+k, ~k=0,1,2,\ldots$.

Regarding details  can be found in $2.2$, 
 Chapter III of  Gelfand and Shilov \cite{Gelfand-Shilov}.

Let $z,w\in\Cx$. $Q^{\alpha~z~w}$ for $\Im z>0$, $\Im w>0$ and $\alpha\in\Cx$ is a distribution defined by analytic continuation from
\bel{Q zv}
Q^{\alpha~z~w}(x,t)~=~\left\{{1\over z |x|^2+ wt^2}\right\}^{\sigma(\alpha)},\qquad \hbox{\small{$\Re\sigma(\alpha)<{n+1\over 2}$}}.
\eeq
From direct computation, we find
\bel{Q ab Transform}
\begin{array}{lr}\ds
\Hat{Q}^{\alpha~\i a~\i b}(\xi,\tau)~=~(-\i)^{\sigma(\alpha)}\iint_{\R^{n+1}} e^{-2\pi\i \big[x\cdot\xi+t\tau\big]} \left\{{1\over a |x|^2+ b t^2}\right\}^{\sigma(\alpha)} dxdt
\\\\ \ds~~
~=~(-\i)^{\sigma(\alpha)}{1\over (\sqrt{a})^n\sqrt{b}}\iint_{\R^{n+1}} e^{-2\pi\i \big[x\cdot\xi/\sqrt{a}+t\tau/\sqrt{b}\big]} \left({1\over  |x|^2+  t^2}\right)^{\sigma(\alpha)} dxdt
\\ \ds~~~~~~~~~~~~~~~~~~~~~~~~~~~~~~~~~~~~~~~~~~~~~~~~~~~~~~~~~~~~~~~~~~~~~
~~~
\hbox{\small{$x\mt x/\sqrt{a}$,~~$t\mt t/\sqrt{b}$}}
\\\\ \ds~~
~=~(-\i)^{\sigma(\alpha)}{\pi^{-{n+1\over 2}+2\sigma(\alpha)}\over (\sqrt{a})^n\sqrt{b}}{\Gamma\left({n+1\over 2}-\sigma(\alpha)\right)\over \Gamma\left(\sigma(\alpha)\right)}\left\{{1\over  |\xi|^2/a+  \tau^2/b}\right\}^{{n+1\over 2}-\sigma(\alpha)},~~~~~\hbox{\small{$0<\Re\sigma(\alpha)<{n+1\over 2}$}}.
\end{array}
\eeq
The last equality is obtained by using 
\bel{power transform}
\int_{\R^\N}e^{-2\pi\i x\cdot\xi} |x|^{\gamma-\N}dx~=~{\pi^{\N-\gamma\over 2}\Gamma\left({\gamma\over 2}\right)\over \pi^{\gamma\over 2}\Gamma\left({\N-\gamma\over 2}\right)} |\xi|^{-\gamma},\qquad 0<\Re\gamma<\N.
\eeq
Replace $a=-\i z$ and $b=-\i w$ inside (\ref{Q ab Transform}). We have
\bel{Q zv Transform}
\begin{array}{lr}\ds
\Hat{Q}^{\alpha~z~w}(\xi,\tau)~=~\iint_{\R^{n+1}} e^{-2\pi\i \left(x\cdot\xi+t\tau\right)} \left\{{1\over z |x|^2+ w t^2}\right\}^{\sigma(\alpha)} dxdt
\\\\ \ds~~~~~~~~~~~~~~~~~~
~=~{\pi^{\sigma(\alpha)-\alpha}\over (\sqrt{ z})^n\sqrt{ w}}{\Gamma\left(\alpha\right)\over \Gamma\left(\sigma(\alpha)\right)}\left\{{1\over  |\xi|^2/z+  \tau^2/w}\right\}^\alpha.
\end{array}
\eeq
\begin{remark}
(\ref{Q zv Transform}) is true  for every $\Im z>0$ and $\Im w>0$.
\end{remark}
For $w=\i b$ fixed, 
both sides of (\ref{Q zv Transform}) are analytic for $\Im z>0$.  Moreover, they are equal at $z=\i a$ for every $a>0$.  Consequently, the two power series of $z$  must have identical coefficients. Given $\Im z>0$, we have a vice versa argument for $w$ where $\Im w>0$. 

On the other hand, we have
\bel{Q ab Transform -}
\begin{array}{lr}\ds
\Hat{Q}^{\alpha~-\i a~-\i b}(\xi,\tau)~=~(-\i)^{-\sigma(\alpha)}{\pi^{-{n+1\over 2}+2\sigma(\alpha)}\over (\sqrt{a})^n\sqrt{b}}{\Gamma\left({n+1\over 2}-\sigma(\alpha)\right)\over \Gamma\left(\sigma(\alpha)\right)}\left\{{1\over  |\xi|^2/a+  \tau^2/b}\right\}^{{n+1\over 2}-\sigma(\alpha)},
\\ \ds~~~~~~~~~~~~~~~~~~~~~~~~~~~~~~~~~~~~~~~~~~~~~~~~~~~~~~~~~~~~~~~~~~~~~~~~~~~~~~~~~~~~~~~~~~~~~~~~~~~~~
\hbox{\small{$0<\Re\sigma(\alpha)<{n+1\over 2}$}}.
\end{array}
\eeq
By replacing $a=\i z$ and $b=\i w$ inside (\ref{Q ab Transform -}),  we find (\ref{Q zv Transform}) again. An analogue to the argument below {\bf Remark B.1} shows (\ref{Q zv Transform}) hold for $\Im z<0, \Im w<0$.

Consider $z=-1+\i a$ and $w=1+\i b$. We have
\bel{Q zv Transform +}
\begin{array}{lr}\ds
\Hat{Q}^{\alpha~z~w}(\xi,\tau)~=~\iint_{\R^{n+1}} e^{-2\pi\i \left(x\cdot\xi+t\tau\right)} \left\{{1\over (-1+\i a) |x|^2+ (1+\i b) t^2}\right\}^{\sigma(\alpha)} dxdt
\\\\ \ds~~~~~~~~~~~~~~~~~~~
~=~{\pi^{\sigma(\alpha)-\alpha}\over (-1+\i a)^{n\over2}(1+\i b)^{1\over 2}}{\Gamma\left(\alpha\right)\over \Gamma\left(\sigma(\alpha)\right)}\left\{{(-1+\i a)(1+\i b)\over  (1+\i b)|\xi|^2+ (-1+\i a) \tau^2}\right\}^\alpha
\end{array}
\eeq
for $0<\Re\sigma(\alpha)<{n+1\over 2}$. Let
\bel{rho a}
\rho(a)~=~\sqrt{1+a^2},\qquad \cos\theta(a)~=~-1/\rho(a),\qquad \sin\theta(a)~=~a/\rho(a).
\eeq
Write $(-1+\i a)^{n\over 2}=\rho(a)^{n\over 2} e^{\i \theta \left({\n\over 2}\right)}$. Clearly, $(-1+\i a)^{n\over 2}\mt e^{\i\pi \left({n\over 2}\right)}$ as $a\mt0$.

On the other hand, choose $z=-1-\i a$ and $w=1-\i b$.  We find
\bel{Q zv Transform -}
\begin{array}{lr}\ds
\Hat{Q}^{\alpha~z~w}(\xi,\tau)~=~\iint_{\R^{n+1}} e^{-2\pi\i \left(x\cdot\xi+t\tau\right)} \left\{{1\over (-1-\i a) |x|^2+ (1-\i b) t^2}\right\}^{\sigma(\alpha)} dxdt
\\\\ \ds~~~~~~~~~~~~~~~~~~~
~=~{\pi^{\sigma(\alpha)-\alpha}\over (-1-\i a)^{n\over2}(1-\i b)^{1\over 2}}{\Gamma\left(\alpha\right)\over \Gamma\left(\sigma(\alpha)\right)}\left\{{(-1-\i a)(1-\i b)\over  (1-\i b)|\xi|^2+ (-1-\i a) \tau^2}\right\}^\alpha
\end{array}
\eeq
for $0<\Re\sigma(\alpha)<{n+1\over 2}$. 
In particular, $(-1-\i a)^{n\over 2}\mt e^{-\i\pi \left({n\over 2}\right)}$ as $a\mt0$.

Consider $\Hat{Q}^{\alpha~z~w}$, $\alpha\in\Cx$  as a distribution defined by analytic continuation from (\ref{Q zv Transform}) where $\Im z\neq0$ and $\Im w\neq0$.

Recall (\ref{P,R,ab})-(\ref{PhiPsi}). 
By taking $a\mt0, b\mt0$, we have
\bel{Phi, R, Psi, P, Transform}
\begin{array}{lr}\ds
\Hat{\Phi}^\alpha~=~\pi^{\sigma(\alpha)-\alpha}e^{-\i\pi\left({ n\over 2}\right)}~{\Gamma\left(\alpha\right)\over \Gamma\left(\sigma(\alpha)\right)}~\Hat{\hbox{R}}^\alpha,
\qquad
\Hat{\Psi}^\alpha~=~\pi^{\sigma(\alpha)-\alpha}e^{\i\pi\left({ n\over 2}\right)}~{\Gamma\left(\alpha\right)\over \Gamma\left(\sigma(\alpha)\right)}~\Hat{\hbox{P}}^\alpha.
\end{array}
\eeq
Recall $\Hat{\hbox{P}}^\alpha, \Hat{\hbox{R}}^\alpha$ and $\Phi^\alpha, \Psi^\alpha$  given at (\ref{PR}) and (\ref{PhiPsi})  in terms of $\Pi^\alpha, \Lambda^\alpha$ and $\hbox{U}^\alpha, \hbox{V}^\alpha$ respectively. 
From direct computation, we find
\bel{Lambda computa}
\begin{array}{lr}\ds
\Hat{\Lambda}^\alpha\left[ e^{i\pi\sigma(\alpha)}-e^{-i\pi\sigma(\alpha)}\right]
~=~\Gamma^{-1}(1-\sigma(\alpha))\left[ e^{\i\pi\sigma(\alpha)} \Hat{\Phi}^\alpha-e^{-\i\pi\sigma(\alpha)}\Hat{\Psi}^\alpha\right]
\\\\ \ds
~=~\Gamma^{-1}(1-\sigma(\alpha)) \Gamma^{-1}(\sigma(\alpha)) \pi^{\sigma(\alpha)-\alpha}\Gamma(\alpha) \left[ e^{-\i\pi\big[{n\over 2}-\sigma(\alpha)\big]} \Hat{\hbox{R}}^\alpha-e^{\i\pi\big[{n\over 2}-\sigma(\alpha)\big]} \Hat{\hbox{R}}^\alpha\right]\qquad\hbox{\small{by (\ref{Phi, R, Psi, P, Transform})}}.
\end{array}
\eeq
Lastly, by using (\ref{Lambda computa}) and the identity $\Gamma(1-z)\Gamma(z)={\pi\over \sin \pi z},~z\in\Cx$, we obtain
\bel{Lambda Transform Formula}
\begin{array}{lr}\ds
\Hat{\Lambda}^\alpha~=~\pi^{\sigma(\alpha)-\alpha-1} \Gamma\left(\alpha\right){1\over 2\i} \Bigg\{- e^{\i\pi \big[{n\over 2}-\sigma(\alpha)\big]} \Hat{\hbox{P}}^\alpha+e^{-\i\pi \big[{n\over 2}-\sigma(\alpha)\big] }\Hat{\hbox{R}}^\alpha\Bigg\}
\\\\ ~~~~~
~=~\pi^{\sigma(\alpha)-\alpha-1}\Gamma(\alpha)\Gamma\left(1-\alpha\right)
\Bigg\{ \Hat{\hbox{U}}^\alpha-\sin\pi\left(\alpha-{1\over 2}\right)\Hat{\hbox{V}}^\alpha\Bigg\}.
\end{array}
\eeq

\small{Department of Mathematics, Westlake University}\\      \small{wangzipeng@westlake.edu.cn}

\end{document}